\documentclass[9pt,amstex]{article}
\usepackage{amssymb,amsmath} 
\usepackage{latexsym}

\newcommand{\ds}{\displaystyle}

\newtheorem{dfn}{Definition}[section] 
\newtheorem{rmk}{Remark}[section]
\newtheorem{thm}{Theorem}[section] 
\newtheorem{cor}{Corollary}[section]
\newtheorem{prop}{Proposition}[section] 
\newtheorem{lem}{Lemma}[section]

\def \als{\mathop{\rightharpoonup}\limits}

\def \ars{\mathop{\leftharpoonup}\limits}
\def\buildrel#1_#2^#3{\mathrel{\mathop{\kern 0pt#1}\limits_{#2}^{#3}}}

\newcommand{\Pf}{{\em Proof}. }
\newcommand{\EPf}
{%
\mbox{}%
\nolinebreak%
\hfill%
\rule{2mm}{2mm}%
\medbreak%
\par%
}

\newcommand{\End}{\mbox{$\mathtt{End}$}}

\newcommand{\id}{\mbox{$\mathtt{Id}$}}

\newcommand{\C}{\mathbb C} 

\newcommand{\R}{\mathbb R}

\newcommand{\g}{{\mathfrak{g}}{}} 
 
\renewcommand{\l}{{\mathfrak{l}}{}}

\newcommand{\p}{{\mathfrak{p}}{}} 
\newcommand{\q}{{\mathfrak{q}}{}}

\newcommand{\s}{{\mathfrak{s}}{}}

\renewcommand{\d}{{\mathfrak{d}}{}} 
 
\renewcommand{\a}{{\mathfrak{a}}{}}

\newcommand{\pw}{{\bf w}{}}

\def\cref#1{Corollary~\ref{#1}}

\title{Universal Deformation Formulae for Three-Dimensional 
Solvable Lie groups}

\author{{\bf Pierre Bieliavsky}\\
Universit\'e Libre de Bruxelles, Belgium,\\
pbiel@ulb.ac.be\\
{\bf Philippe Bonneau}\\
Universit\'e de Bourgogne, France,\\
bonneau@u-bourgogne.fr\\
{\bf Yoshiaki Maeda}\\
Keio University, Japan,\\
maeda@math.keio.ac.jp}
\date{June 10, 2003} 

\begin{document}
\maketitle

\begin{abstract}
We apply methods from strict quantization of solvable
symmetric spaces to obtain universal deformation formulae for actions
of every three-dimensional solvable Lie group.  We also study compatible 
co-products by generalizing the notion of smash product in the context
of Hopf algebras. We investigate in particular 
the dressing action of the `book' group on $SU(2)$.
\end{abstract}

\section{Introduction}
Let $G$ be a group acting on a set $M$. Denote by $\tau:G\times M\to 
M:(g,x)\mapsto\tau_g(x)$ the (left) action and by
$\alpha:G\times\mbox{Fun}(M)\to\mbox{Fun}(M)$ the corresponding action
on the space of (complex valued) functions (or formal series) on $M$
($\alpha_g:=\tau^\star_{g^{-1}}$).  Assume that on a subspace
${\mathbb A}\subset\mbox{Fun}(G)$, one has an associative $\C$-algebra
product $\star^G_{\mathbb A}:{\mathbb A}\times {\mathbb
A}\to {\mathbb A}$ such that
\begin{enumerate}
\item[(i)] ${\mathbb A}$ is invariant under the (left) regular action of $G$ on 
$\mbox{Fun}(G)$,
\item[(ii)] the product $\star^G_{\mathbb A}$ is left-invariant as well
i.e. for all $g\in G; a,b\in {\mathbb A}$, one has
\begin{equation*}
(L_g^\star a) \star^G_{\mathbb A} (L_g^\star b)=L_g^\star(a\star^G_{\mathbb A}b).
\end{equation*}
\end{enumerate}

\noindent Given a function on $M$, $u\in\mbox{Fun}(M)$, and a point $x\in M$, one 
denotes by $\alpha^x(u)\in\mbox{Fun}(G)$ the function on $G$ defined as 
\begin{equation*}
\alpha^x(u)(g):=\alpha_g(u)(x).
\end{equation*}
Then one readily observes that the subspace ${\mathbb B}\subset\mbox{Fun}(M)$ defined as 
\begin{equation*}
{\mathbb B}:=\{u\in\mbox{Fun}(M)\, |\, \forall x\in M:\alpha^x(u)\in
{\mathbb A}\}
\end{equation*}
becomes an associative $\C$-algebra when endowed
with the product $\star^M_{\mathbb B}$ given by
\begin{equation}\label{UDF}
u\star^M_{\mathbb B}v(x):=(\alpha^x(u)\star^G_{\mathbb A}\alpha^x(v))(e)
\end{equation}
($e$ denotes the neutral element of $G$).
Of course, all this can be defined for right actions as well. 
\begin{dfn}
Such a pair $({\mathbb A},\star^G_{\mathbb A})$ is called a (left) {\bf universal deformation}
of $G$, while Formula (\ref{UDF}) is called the associated {\bf universal 
deformation formula} (briefly {\bf UDF}).
\end{dfn}
In the present article, we will be concerned with the case where $G$ is a Lie 
group. The function space ${\mathbb A}$ will be either 

\noindent - a functional subspace (or a topological completion) of 
$C^\infty(G,\C)$ containing the smooth compactly supported functions in which case we will talk
about {\bf strict deformation} (following Rieffel \cite{Rieffel1}), 

\noindent or, 

\noindent - the space ${\mathbb A}=C^\infty(G)[[\hbar]]$ of formal power series with coefficients 
in the smooth functions on $G$ in which case, we'll speak about {\bf 
formal deformation}. In any case, we'll assume the product $\star^G_{\mathbb A}$ 
admits an asymptotic expansion of  star-product type:
\begin{equation*}
a\star^G_{\mathbb A}b\sim
ab+\frac{\hbar}{2i}\pw({\rm d}u,{\rm d}v)+o(\hbar^2)\qquad (a,b\in C^\infty_c(G)),
\end{equation*}
where $\pw$ denotes some (left-invariant) Poisson bivector on
$G$\cite{BFFLS}.  In the strict cases considered here, the product
will be defined by an integral three-point kernel $K\in
C^\infty(G\times G\times G)$:
\begin{equation*}
a\star^G_{\mathbb A}b(g):=\int_{G\times G}a(g_1)\,b(g_2)K(g,g_1,g_2){\rm d}g_1\,{\rm 
d}g_2\qquad(a,b\in {\mathbb A})
\end{equation*}
where ${\rm d}g$ denotes a normalized left-invariant Haar 
measure on $G$. Moreover, our kernels will be
of  {\bf WKB type} \cite{Wein, Kara} i.e.:
\begin{equation*}
K=A\,e^{\frac{i}{\hbar}\Phi},
\end{equation*}
with $A$ (the {\bf amplitude}) and $\Phi$ (the {\bf phase}) in $C^\infty(G\times 
G\times G,\R)$ being invariant under the (diagonal) action by 
left-translations.

\noindent Note that in the case where the group $G$  acts smoothly on a smooth
manifold $M$ by diffeomorphisms: $\tau:G\times M\to 
M:(g,x)\mapsto\tau_g(x)$, the first-order expansion term of
$u\star^M_{\mathbb B}v,\quad u,v\in C^{\infty}(M)$ defines a Poisson
structure $\pw^{M}$ on $M$ which can be expressed in terms of a basis $\{X_i\}$ of 
the Lie  algebra $\g$ of $G$ as:
\begin{equation}\label{PS}
    \pw^{M}=\left[\pw_e\right]^{ij}X^\star_i\wedge
X^\star_j,
\end{equation}
where $X^\star$ denotes the fundamental vector field on $M$
associated to $X\in\g$.

\vspace{5mm}

\noindent Strict deformation theory in the WKB context was initiated 
by Rieffel in \cite{Rieffel2} in the cases where $G$ is either Abelian or
1-step nilpotent.  Rieffel's work has led to what is now called 
`Rieffel's machinery'; producing a whole class of exciting non-commutative
manifolds (in Connes sense) from the data of Abelian group actions on
$C^{\star}$-algebras \cite{CoLa}.

\noindent The study of formal UDF's for non-Abelian group actions in
our context was initiated in \cite{GiaZha}
where the case of the group of  affine transformations of the real line
(`$ax+b$') was explicitly described.

\noindent In the strict  (non-formal) setting, UDF's for Iwasawa subgroups
of $SU(1,n)$ have been explicitly given in \cite{BiMas}.  These were
obtained by adapting a method developed by one of us in the symmetric space
framework \cite{Biep00a}. In particular when $n=1$ one obtains strict UDF's
for the group '$ax+b$' (we recall this in the appendix).

\noindent In the present work, we provide in the Hopf algebraic context (strict and formal)
UDF's for every solvable three-dimensional Lie group endowed with any left-invariant Poisson 
structure. The method used to obtain those is based on strict quantization of solvable 
symmetric spaces, on existing UDF's for '$ax+b$' and on a generalization of 
the classical definition of smash products
in the context of Hopf algebras.
As an application we analyze the particular case of the dressing
action of the Poisson dual Lie group of $SU(2)$ when endowed with the
Lu-Weinstein Poisson structure.  

\section{UDF's for three-dimensional solvable Lie groups}
In this section $G$ denotes a three-dimensional solvable Lie group
with Lie algebra $\g$.
\begin{dfn}
Let $\pw$ be a left-invariant Poisson bivector field on $G$. The pair 
$(G,\pw)$ is called a {\bf pre-symplectic} Lie group.
\end{dfn}
The terminology 'symplectic' is after Lichn\'erowicz \cite{Li} who studied
this type of structures.

\noindent One then observes
\begin{prop}
\begin{enumerate}
\item[(i) ]The orthodual $\s$ of the radical of $\pw_e$
\begin{equation*}
\s:=(\mbox{rad }\pw_e)^\perp
\end{equation*}
is a subalgebra of $\g$.
\item[(ii)] Every two-dimensional subalgebra of $\g$ can be seen as the orthodual 
of the radical of a Poisson bivector.
\item[(iii)] The symplectic leaves of $\pw$ are the left classes of the 
analytic subgroup $S$ of $G$ whose Lie algebra is $\s$.
\end{enumerate}
\end{prop}
We now fix a pair $(\g,\pw_e)$ to be such a pre-symplectic Lie algebra with 
associated symplectic Lie algebra $(\s,\pw_e|_{\s^\star\times\s^\star})$.
We assume 
\begin{equation*}
\dim\s=2.
\end{equation*}
We denote by $\d$ the first derivative of $\g$:
\begin{equation*}
\d:=[\g,\g].
\end{equation*}
Note that $\d$ is an Abelian algebra.

\subsection{Case 1: $\dim \d=2$}
In this case $\g$ can be realized as a split extension of Abelian algebras:
$$
0\to\d\longrightarrow\g{\longrightarrow}\a\to0,
$$
with $\dim\a=1$. We denote by 
\begin{equation*}
\rho:\a\to\End(\d)
\end{equation*}
the splitting homomorphism.
\subsubsection{$\s\cap\d\neq\d$}\label{SNEQD}
In this case one sets
\begin{equation*}
\p:=\s\cap\d,
\end{equation*}
and one may assume
\begin{equation*}
\a\subset\s.
\end{equation*}
Note that $\s$ is then isomorphic to the Lie algebra of the group '$ax+b$'.
\subsubsection*{\ref{SNEQD}.1 The representation $\rho$ is semisimple}
In this case one has $\q\subset\d$ such that
\begin{eqnarray*}
\q\oplus\p=\d\mbox{ and}\\
\left[ \a,\q \right] \subset\q.
\end{eqnarray*}
Let $Q$ and $S$ denote the analytic subgroups of $G$ whose Lie algebras 
are $\q$ and $\s$ respectively. Consider the mapping
\begin{equation}\label{CHART}
Q\times S\rightarrow G:(q,s)\mapsto qs.
\end{equation}
Assume this is a global diffeomorphism.  Endow the symplectic Lie group $S$
with a left-invariant deformation quantization $\star^T$ as constructed in 
\cite{BiMas} or \cite{BiMae} (see also the appendix of the present paper).
One then readily verifies
\begin{thm}\label{INDUC}
Let $({\bf H}^S,\star^T)$ be an associative algebra of functions on $S$ such 
that ${\bf H}^S\subset \mbox{Fun}(S)$ is an invariant subspace w.r.t. the 
left regular representation of $S$ on $\mbox{Fun}(S)$ and $\star^T$ is a
$S$-left-invariant product on ${\bf H}^S$. Then the function space 
${\bf H}:=\{u\in\mbox{Fun}(G)\ |\ \forall p\in Q:u(q,.)\in {\bf H}^S\}$ is an 
invariant subspace of $\mbox{Fun}(G)$ w.r.t. the left regular representation 
of $G$ and the formula
\begin{equation*}
u\star v(q,s):=(u(q,.)\star^Tv(q,.))(s)
\end{equation*}
defines a $G$-left-invariant associative product on ${\bf H}$. 
\end{thm}

\subsubsection*{\ref{SNEQD}.2 The representation $\rho$ is not semisimple}
Since it cannot be nilpotent one can find bases $\{A\}$ of $\a$
and $\{X,Y\}$ of $\d$ such that
\begin{equation*}
\mbox{span}\{X\}=\p
\end{equation*}
and 
\begin{equation*}
\rho(A)=\left(
\begin{array}{cc}
1&\lambda\\
0&1
\end{array}
\right)\qquad \lambda\in\R_0.
\end{equation*}
\begin{lem}
$\g$ can be realized as a subalgebra of the transvection algebra of a four-dimensional
symplectic symmetric space $M$.
\end{lem}
\Pf
The symplectic triple associated to the symplectic symmetric space (cf. \cite{Biep95a} or \cite{Biep00a} Definition 2.5) is 
$(\mathfrak{G},\sigma,\Omega)$ where 
\begin{eqnarray*}
\mathfrak{G}&=&\mathfrak{K}\oplus\mathfrak{P}\\
\mathfrak{K}&=&\mbox{span}\{k_1,k_2\}\\
\mathfrak{P}&=&\mbox{span}\{e_1,f_1,e_2,f_2\};
\end{eqnarray*}
with Lie algebra table given by
\begin{eqnarray*}
\left[f_1,k_1\right]&=&e_1\\
\left[f_1,k_2\right]&=&e_2+\lambda e_1\\
\left[f_1,e_1\right]&=&k_1\\
\left[f_1,e_2\right]&=&k_2+\lambda k_1;
\end{eqnarray*}
$f_2$ being central.
The symplectic structure on $\mathfrak{P}$ is given by
\begin{equation*}
\Omega=\alpha e_1^*\wedge f_1^*+\beta e_2^*\wedge f_2^*\qquad(\alpha,\beta\neq0).
\end{equation*}
The subalgebra 
\begin{equation*}
\mbox{span}\{X:=\frac{1}{2}(k_1+e_1),\,Y:=\frac{1}{2}(k_2+e_2),\,A:=f_1,\}
\end{equation*}
is then isomorphic to $\g$.
\EPf

\noindent Note that the transvection subalgebra
\begin{equation*}
\mathfrak{S}=\g\oplus\R f_2
\end{equation*}
is transverse to the isotropy algebra $\mathfrak{K}$ in $\mathfrak{G}$.
As a consequence the associated connected simply connected Lie group
$\bf S$ acts acts transitively on the symmetric space $M$. Since 
$\dim M=\dim{\bf S}=4$, one may identify $M$ with the group manifold of $\bf S$.
It turns out that the symplectic symmetric space $M$ admits a strict deformation
quantization which is transvection invariant \cite{Biep00a}. Therefore one has
\begin{thm}
The symplectic Lie group $\bf S$ admits a left-invariant strict deformation
quantization. Since $G$ is a group direct factor of $\bf S$, the deformed 
product on $\bf S$ restricts to $G$ providing a UDF for $G$.
\end{thm}

\subsubsection{$\s=\d$}\label{ABELIAN}
In this case, a UDF is obtained by writing Moyal's formula where one replaces 
the partial derivatives $\partial_i$'s by (commuting) left-invariant vector fields
$\widetilde{X_i}$ associated with generators $X_i$ of $\d$:
for $a,b\in C^\infty(G)[[\hbar]]$, one sets
\begin{equation*}
a\star_\hbar b:=a.\exp\left( \pw_e|_\d^{ij}\stackrel{\leftarrow}{\widetilde{X_i}}\wedge \stackrel{\rightarrow}{\widetilde{X_j}}\right).b
\end{equation*}

\subsection{Case 2: $\dim\d=1$}
In this case $\g$ is isomorphic to either the Heisenberg algebra or to the direct sum 
of $\R$ with the Lie algebra of '$ax+b$'.  The first case has been studied by Rieffel in \cite{Rieffel1}.
Let then 
\begin{equation*}
\g=\mathfrak{L}\oplus\R Z
\end{equation*}
where $\mathfrak{L}=\mbox{span}\{A,E\}$ with $[A,E]=E$. Up to automorphism one has either 
$\s=\mathfrak{L}$, $\s=\mbox{span}\{E,Z\}$ or $\s=\mbox{span}\{A,Z\}$. The last ones are Abelian
thus UDF's in these cases are obtained the same way as in Subsection \ref{ABELIAN}.
The first case reduces to an action of '$ax+b$' and has therefore been treated in \cite{BiMas}
or \cite{BiMae}.

\section{Crossed, Smash and Co-products}
\indent Every algebra, coalgebra, bialgebra, Hopf algebra and vector space is 
taken over the field $k=\mathbb{R}$ or $\mathbb{C}$.  For classical
definitions and facts on these subjects, we refer to \cite{Swem69a}, 
\cite{Abee80a} or more fundamentally to \cite{MiMo65}.

\noindent To calculate with a coproduct $\Delta$, we use the Sweedler notation
\cite{Swem69a}: 
$\ds
\Delta(b)=\sum_{(b)} b_{(1)}\otimes b_{(2)}.
$ 

For the classical definitions of a $B$-module algebra,
coalgebra and bialgebra we refer to \cite{Abee80a}.

\medskip

The literature on Hopf algebras contains a large collection of what we
can call generically semi-direct products, or crossed products. Let us describe 
some of them. The simplest example of these crossed products is usually called 
the smash product (see \cite{Swem68a, Mol77a}):

\begin{dfn}\label{smash}
Let $B$ be a bialgebra and $C$ a $B$-module algebra. The smash product
$C\sharp B$ is the algebra constructed on the vector space $C\otimes B$ 
where the multiplication  is defined as 
\begin{equation}
(f\otimes a)\stackrel{\rightharpoonup}{\star}(g\otimes b)
=\sum_{(a)}f\ (a_{(1)}\rightharpoonup g)\otimes a_{(2)} b\ ,\qquad
f,g\in C\ ,\quad a,b\in B.
\end{equation}

\end{dfn}

Assuming $B$ is cocommutative, we now introduce a generalization of the 
smash product.
\begin{dfn}\label{LRsmash}
Let $B$ be a cocommutative bialgebra and $C$ a $B$-bimodule algebra 
(i.e. a $B$-module algebra for both, left and right, $B$-module structures). 
The {\bf L-R-smash product}
$C\natural B$ is the algebra constructed on the vector space $C\otimes B$ 
where the multiplication is defined by 
\begin{equation}
(f\otimes a)\star (g\otimes b)=\sum_{(a)(b)}
                (f\leftharpoonup b_{(1)}) (a_{(1)}\rightharpoonup g)
                                        \otimes a_{(2)} b_{(2)}\ ,\qquad
f,g\in C\ ,\quad a,b\in B.
\end{equation}

\end{dfn}

\begin{rmk} 
{\rm The fact that the L-R-smash product $C\natural B$ is an associative algebra 
follows by an easy adaptation of the proof of the associativity
for the smash product \cite{Mol77a}}.
\end{rmk} 

In the same spirit, one has
\begin{lem}\label{cop}
If $C$ is a $B$-bimodule bialgebra, the natural tensor 
product coalgebra structure on $C\otimes B$ defines a bialgebra 
structure to $C\natural B$. 

If $C$ and $B$ are Hopf algebras, $C\natural B$ is a Hopf
algebra as well, defining the antipode by
\begin{eqnarray}
J_\star(f\otimes a) & = & \sum_{(a)}J_B(a_{(1)})\als J_C(f)\ars J_B(a_{(2)})
	\otimes J_B(a_{(3)})\\ 
	& = & \sum_{(a)}(1_C\otimes J_B(a_{(1)}))\star (J_C(f)\otimes 1_B)\star 
	(1_C\otimes J_B(a_{(2)})).\nonumber
\end{eqnarray}
\end{lem} 

\noindent Now by a careful computation, one proves
\begin{prop}\label{*S}
Let $B$ be a cocommutative bialgebra, $C$ be a $B$-bimodule algebra
and $(C\natural B, \star )$ be their L-R-smash product.\\ 
Let $T$ be a linear automorphism of $C$ (as vector space). We define:
\begin{enumerate}
\item[(i)]  the product $\bullet^T$ on $C$ by 
\begin{equation*}
f\bullet^T g=T^{-1}\left(T(f).T(g)\right);
\end{equation*}
\item[(ii)]  the left and right $B$-module structures, $\als^T$ and 
$\ars^T$, by 
\begin{equation*}
a\als^T f := T^{-1}\left(a\rightharpoonup T(f)\right)\ \mbox{ and }\
f\ars^T a := T^{-1}\left(T(f)\leftharpoonup a\right);
\end{equation*}
\item[(iii)]  the product, $\star^T$, on $C\otimes B$ by 
\begin{equation*}
(f\otimes a)\star^T (g\otimes b)
={\bf T}^{-1}\left( {\bf T}(f\otimes a) \star {\bf T}(g\otimes b)\right)
\qquad\mbox{ where } {\bf T}:=T\otimes \id .
\end{equation*}

\end{enumerate} 
Then $(C, \bullet^T)$ is a $B$-bimodule algebra for $\als^T$ and $\ars^T$ 
and $\star^T$ is the L-R-smash product defined by these structures.

\noindent Moreover, 
if $(C, . , \Delta_C , J_C , \rightharpoonup , \leftharpoonup )$ 
is a Hopf algebra  and a $B$-module bialgebra, then 
$$C_T:=(C, \bullet^T , \Delta^T_C:=(T^{-1}\otimes T^{-1})\circ\Delta_C\circ T ,
J^T_C:=T^{-1}\circ J_C\circ T , \als^T , \ars^T )$$ 
is also a Hopf algebra  and a $B$-module bialgebra. 
Therefore, by Lemma \ref{cop}, 
$$(C_T\natural B , \star^T , \Delta^T = (23)\circ (\Delta^T_C\otimes 
\Delta_B) , J_\star^T ),$$
is a Hopf algebra, 
$\Delta^T$ being the natural tensor product coalgebra structure
on $C_T\natural B$ (with $(23):C\otimes C\otimes B\otimes B\rightarrow 
C\otimes B\otimes C\otimes B,\
c_1\otimes c_2\otimes b_1\otimes b_2\mapsto
c_1\otimes b_1\otimes c_2\otimes b_2$)
and $J_\star^T$ being the antipode given on $C_T\natural B$ by Lemma \ref{cop}.
Also, one has $$\Delta^T=({\bf T}^{-1}\otimes {\bf T}^{-1})\circ (23)
\circ (\Delta_C\otimes\Delta_B)\circ {\bf T}\ \mbox{ and }\
J_\star^T={\bf T}^{-1}\circ J_\star \circ {\bf T}\qquad \mbox{with } {\bf T}=T\otimes \id .$$ 
\end{prop} 

\section{Examples in deformation quantization}\label{EX}

\subsection{A construction on $T^\star(G)$}
Let  $G$ be a Lie group with Lie algebra $\g$ and $T^\star(G)$ its cotangent 
bundle. We denote by $\mathsf{U}\mathfrak{g}$,
$\mathsf{T}\mathfrak{g}$ and $\mathsf{S}\mathfrak{g}$ respectively the 
enveloping, tensor and symmetric algebras of $\mathfrak{g}$. 
Let $\mathsf{Pol}(\mathfrak{g}^*)$ be
the algebra of polynomial functions on $\mathfrak{g}^*$. 
We have the usual identifications:
$$\mathcal{C}^\infty (T^*G)\simeq \mathcal{C}^\infty (G\times \mathfrak{g}^*)
\simeq \mathcal{C}^\infty(G)\hat{\otimes} 
\mathcal{C}^\infty(\mathfrak{g}^*) \supset 
\mathcal{C}^\infty(G)\otimes\mathsf{Pol}(\mathfrak{g}^*)
\simeq \mathcal{C}^\infty(G)\otimes\mathsf{S}\mathfrak{g}.$$

\noindent First we deform $\mathsf{S}\mathfrak{g}$ via 
the ``parametrized version'' $\mathsf{U}_t\mathfrak{g}$ 
of $\mathsf{U}\mathfrak{g}$ defined by
$$\mathsf{U}_t\mathfrak{g}= \frac{\mathsf{T}\mathfrak{g}[[t]]}
{<XY-YX-t[X,Y]; X,Y\in \mathfrak{g}>}.$$
$\mathsf{U}_t\mathfrak{g}$ is naturally a Hopf algebra with 
$\Delta(X)=1\otimes X + X\otimes 1$, $\epsilon(X)=0$ and 
$S(X)=-X$ for $X\in \mathfrak{g}$.
For $X\in\g$, we denote by $\widetilde{X}$ (resp. $\overline{X}$) the
left- (resp. right-) invariant vector field on $G$ such that 
$\widetilde{X}_e=\overline{X}_e=X$. We consider
the following $k[[t]]$-bilinear actions of $B=\mathsf{U}_t\mathfrak{g}$
on $C=\mathcal{C}^\infty(G)[[t]]$, for $f\in C$ and $\lambda\in [0,1]$:
\begin{enumerate}
\item[(i)] $(X \rightharpoonup f)(x) = t (\lambda -1)\ 
(\widetilde{X}.f)(x)$,
\item[(ii)] $(f \leftharpoonup X)(x) = t \lambda \ (\overline{X}.f)(x)$.
\end{enumerate}
One then has
\begin{lem}
$C$ is a $B$-bimodule algebra 
w.r.t. the above left and right actions (i) and (ii).
\end{lem}

\begin{dfn}
We denote by $\star_\lambda$ the star product on
$\left(\mathcal{C}^\infty(G)\otimes\mathsf{Pol}(\mathfrak{g}^*)\right)[[t]]$ 
given by the L-R-smash product on 
$\mathcal{C}^\infty(G)[[t]]\otimes\mathsf{U}_t\mathfrak{g}$
constructed from the bimodule structure of the preceding lemma.
\end{dfn}

\begin{prop}
For $G=\mathbb{R}^n$, $\star_{\frac{1}{2}}$ is the Moyal star product 
(Weyl ordered), $\star_{0}$ is the standard ordered star product and 
$\star_{1}$ the anti-standard ordered one. 
In general $\star_\lambda$ yields the 
$\lambda$-ordered quantization, within the 
notation of M. Pflaum \cite{Pflm99a}. 
\end{prop}

\begin{rmk}
{\rm 
In the general case, it would be interesting to 
compare our $\lambda$-ordered L-R smash product with classical
constructions of star products on $T^\star(G)$, with Gutt's product
as one example \cite{Guts83a}.
}
\end{rmk}

\subsection{Hopf structures}\label{CP}
We have discussed (see Lemma \ref{cop})
the possibility of having a Hopf structure on $C\natural B$.  Let
us consider the particular case of $\mathcal{C}^\infty
(\mathbb{R}^n)[[t]] \natural\mathsf{U}_t\mathbb{R}^n =
\mathcal{C}^\infty (\mathbb{R}^n)[[t]] \natural\mathsf{S}\mathbb{R}^n$
($\mathbb{R}^n$ is commutative).  $\mathsf{S}\mathbb{R}^n$ is endowed
with its natural Hopf structure but we also need a Hopf
structure on $\mathcal{C}^\infty (\mathbb{R}^n)[[t]] =
\mathcal{C}^\infty (\mathbb{R}^n)\otimes\mathbb{R}[[t]].$ We will not
use the usual one. Our alternative structure is defined as follows.
\begin{dfn}
We endow $\mathbb{R}[[t]]$ with the usual product, the co-product
$\Delta(P)(t_1,t_2):=P(t_1+t_2)$, the co-unit $\epsilon(P)=P(0)$
and the antipode $J(t)=-t$.
We consider the Hopf algebra
$(\mathcal{C}^\infty (\mathbb{R}^n), . , 
{\bf 1}, \Delta_C , \epsilon_C , J_C )$, with 
pointwise multiplication, the unit ${\bf 1}$ 
(the constant function of value $1$), the coproduct
$\Delta_C (f)(x,y)=f(x+y)$, the co-unit $\epsilon (f)=f(0)$
and the antipode $J_C(f)(x)=f(-x)$. 
The tensor product of these two Hopf algebras
then yields a  Hopf algebra denoted by
$$(\mathcal{C}^\infty (\mathbb{R}^n)[[t]], . , {\bf 1} ,
\Delta_t , \epsilon_t , J_t).$$ 
Note that $\Delta_t$ and $J_t$ are not linear in $t$.
We then define, on the L-R smash $\mathcal{C}^\infty (\mathbb{R}^n)[[t]]
\natural\mathsf{S}\mathbb{R}^n$,
$$\Delta_\star := (23)\circ (\Delta_t\otimes\Delta_B),\quad
\epsilon_\star := \epsilon_t\otimes\epsilon_B
\quad \mbox{ and }\quad J_\star\ \mbox{ as in Lemma \ref{cop} }.$$ 
\end{dfn}

\begin{prop}
($\mathcal{C}^\infty (\mathbb{R}^n)[[t]]\natural\mathsf{S}\mathbb{R}^n,
\star_\lambda , {\bf 1}\otimes 1 , \Delta_\star , \epsilon_\star , J_\star )$
is a Hopf algebra.
\end{prop}

\begin{rmk}
{\rm 
The case $\lambda=\frac{1}{2}$ yields the usual Hopf structure on 
the enveloping algebra of the Heisenberg Lie algebra.
}
\end{rmk}

\section{The `Book' Algebra}
\begin{dfn}
The {\bf book Lie algebra} is the three-dimensional Lie algebra
$\g$ defined as the split extension of Abelian Lie algebras $\a:=\R.A$
and $\d:=\R^2$:
\begin{equation*}
0\to\d\to\g\to\a\to0,
\end{equation*}
where the splitting homomorphism $\rho:\a\to\End(\d)$ is defined by 
\begin{equation*}
\rho(A):=\id_\d.
\end{equation*}
\end{dfn}
The name comes from the fact that the regular co-adjoint orbits in $\g^\star$
are open half planes sitting in $\g^\star=\R^3$ ressembling the 
pages of an open book.

\noindent We are particularly interested in this example because the 
associated connected simply connected Lie group $G$
turns out to be the solvable Lie group underlying the Poisson dual
of $SU(2)$ when endowed with the Lu-Weinstein Lie-Poisson structure \cite{LW}.
Explicitly, one has the following situation (see \cite{LuRat} or \cite{Slee}).  Set
$\mathfrak{K}:=\mathfrak{su}(2)$ and consider the Cartan decomposition
of $\mathfrak{K}^{\C}=\mathfrak{sl}_2(\C)$:
\begin{equation*}
\mathfrak{K}^{\C}=\mathfrak{K}\oplus i\mathfrak{K}.
\end{equation*}
Fix a maximal Abelian subalgebra $\a$ in $i\mathfrak{K}$, consider the 
corresponding root space decomposition
\begin{equation*}
\mathfrak{K}^{\C}=\mathfrak{K}^+\oplus\a^\C\oplus\mathfrak{K}^-,
\end{equation*}
and set
\begin{equation*}
\d:=\mathfrak{K}^+\simeq\C\simeq\R^2.
\end{equation*}
The Iwasawa part $\a\oplus\mathfrak{K}^+$ is then a realization of the book 
algebra $\g=\a\times\d$. One denotes by
\begin{equation*}
K^\C=KG
\end{equation*}
the Iwasawa decomposition of $K^\C$. For an element $\gamma\in K^\C$, 
we have the corresponding factorization:
\begin{equation*}
\gamma=\gamma_K\gamma_G.
\end{equation*}
The dressing action
\begin{equation}\label{DRESS}
G\times K\stackrel{\tau}{\longrightarrow}K 
\end{equation}
is then given by $$\tau_g(k):=(gk)_K.$$ 
In the same way, denote by
\begin{equation*}
\mathfrak{K}^\C\to\mathfrak{K}:A\mapsto A_{\mathfrak{K}}
\end{equation*}
the projection parallel to $\g$. The infinitesimal 
dressing action of $G$ on $K$ is now given by
\begin{equation*}
X^\star_k=-L_{k_{\star_e}}\left[\left(\mbox{Ad}(k^{-1})X\right)_\mathfrak{K}\right],\quad X\in\g.
\end{equation*}
A choice of a pre-symplectic structure on $G$---or equivalently (up to a 
scalar) a choice of a two-dimensional Lie algebra $\s$ in $\g$--- then 
yields, 
via the dressing action (\ref{DRESS}), a Poisson structure $\pw^\s$ on $K$
(cf.  (\ref{PS})):
\begin{equation*}
L_{k^{-1}_{\star_k}}(\pw^\s_k):=
\left(\mbox{Ad}(k^{-1})S_1\right)_\mathfrak{K}\wedge
\left(\mbox{Ad}(k^{-1})S_2\right)_\mathfrak{K}
\end{equation*}
where $\{S_{1},S_2\}$ is a basis of $\s$. Note that this Poisson structure 
is compatible with the Lu-Weinstein structure.

\noindent For each choice of symplectic Lie subalgebra $\s$ of $\g$, the preceding section then produces deformations of
$(K,\pw^\s)$.  

\section*{Appendix: UDF's for '$ax+b$'}
Let $S='ax+b'$ denote the two-dimensional solvable Lie group presented as follows.
As a manifold, one has
\begin{equation*}
S=\R^2=\{(a,\ell)\},
\end{equation*}
and the multiplication law is given by
\begin{equation*}
(a,\ell).(a',\ell'):=(a+a',\,e^{-a'}\ell+\ell').
\end{equation*}
Observe that the symplectic form on $S$ defined by
\begin{equation*}
\omega:=da\wedge dl
\end{equation*}
is left-invariant.

\noindent We now define two specific real-valued functions.
First, the two-point function $A\in C^\infty(S\times S)$ given by
\begin{equation*}
A(x_1,x_2):=\cosh(a_1-a_2)
\end{equation*}
where $x_i=(a_i,\ell_i)\quad i=1,2$.
Second, the three-point function $\Phi\in C^\infty(S\times S\times S)$ given by
\begin{equation*}
\Phi(x_0,x_1,x_2):=\oint_{0,1,2}\sinh(a_0-a_1)\ell_2
\end{equation*}
where $\oint$ stands for cyclic summation. One then has
\begin{thm}\cite{Biep00a,BiMas}
There exists a family of functional subspaces 
\begin{equation*}
\{\cal H\}_{\hbar\in\R}\subset C^\infty(S)
\end{equation*}
such that
\begin{enumerate}
\item[(i)] for all $\hbar\in\R$,
$$C^\infty_c(S)\subset{\cal H}_\hbar;$$
\item[(ii)] For all $\hbar\in\R\backslash\{0\}$ and  $u,v\in C^\infty_c(M)$, the formula:
\begin{equation}\label{PROD}
u\star_\hbar v(x_0):=\int_{M\times M}u(x_1)
\,v(x_2)\,A(x_1,x_2)\,e^{\frac{i}{\hbar}\Phi(x_0,x_1,x_2)}\,{\rm d}x_1\,{\rm d}x_2
\end{equation}
extends as an associative product on ${\cal H}_\hbar$ ($\,{\rm d}x$ 
denotes some normalization of the 
symplectic volume on $(S,\omega)$). Moreover, (for 
suitable $u,v$ and $x_{0}$) a stationary phase method yields a power 
series expansion of the form
\begin{equation}\label{STAR}
u\star_\hbar v(x_0)\sim 
uv(x_0)+\frac{\hbar}{2i}\{u,v\}(x_{0})+o(\hbar^{2});
\end{equation}
where $\{\,,\,\}$ denotes the symplectic Poisson bracket on $(S,\omega)$.
\item[(iii)] The pair $({\cal H}_\hbar,\star_\hbar)$ is a topological (pre) Hilbert 
algebra on which the group $S$ acts on the left  by 
automorphisms.
\end{enumerate}
\end{thm}

\noindent Now setting $\R^2=\a\times\l;\quad a\in\a,\ell\in\l$, one gets the linear isomorphisms
$$
\mathcal{C}^\infty (S)
\simeq \mathcal{C}^\infty (\mathfrak{l}) \hat\otimes
\mathcal{C}^\infty (\mathfrak{a}) \simeq \mathcal{C}^\infty
(\mathfrak{l}) \hat\otimes \mathcal{C}^\infty (\mathfrak{l}^\ast)
\supset\mathcal{C}^\infty
(\mathfrak{l})\otimes Pol(\mathfrak{l}^\ast)
\simeq
\mathcal{C}^\infty (\mathfrak{l})\otimes \mathsf{U}\mathfrak{l}.
$$
The quantization (\ref{PROD}) on such a space turns out to be a L-R-smash 
product. Namely, one has
\begin{prop}
The formal version of the invariant quantization (\ref{PROD})  is 
a L-R-smash product of the form
$\star^T$ (cf. Proposition \ref{*S}).
\end{prop}
\Pf
Let ${\cal S}(\l)$ denote the Schwartz space on the vector space $\l$.
In \cite{Biep00a}, one shows that the map
$$
\begin{array}{ccc}
{\cal S}(\l) &\stackrel{T}{\longrightarrow}&{\cal S}(\l) \\
u&\mapsto&T(u):=F^{-1}\phi_\hbar^\star F(u),
\end{array}
$$
where 
\begin{equation*}
\phi_\hbar:\l^*\to\l^*:a\mapsto\frac{2}{\hbar}\sinh(\frac{\hbar}{2}a),
\end{equation*}
is a linear injection for all $\hbar\in\R$ ($F$ denotes the partial 
Fourier transform  w.r.t. the variable $\ell$).
An asymptotic expansion in a power series in $\hbar$ then yields
a formal equivalence again denoted by ${T}$:
$$
T:=\id+o(\hbar):C^\infty(\l)[[\hbar]]\to C^\infty(\l)[[\hbar]].
$$
Carrying the Moyal star product on $(\a\times\l,\omega)$ by
$\ds \ {\bf T}:=T\otimes \id \ $ yields a star product on $S$ 
which coincides with
the asymptotic expansion of the left-invariant star-product (\ref{PROD}) \cite{Biep00a,BiMas}.
\EPf
\noindent Subsection \ref{CP} then yields
\begin{cor}
The (formal) UDF (\ref{PROD}) admits  compatible co-product and antipode.
\end{cor}

\end{document}